\documentclass[a4paper,12pt,draft]{article}
\usepackage[cp1251]{inputenc}
\usepackage[english,russian]{babel}
\usepackage{amscd,amssymb,amsfonts,amsmath,array}

\newcommand{\rx}{{\rm x}}

\newcommand{\ru}{{\rm u}}
\newcommand{\rz}{{\rm z}}

\begin{document}  

\title{\bf КРИТЕРИЙ	ВНУТРЕННЕЙ УСТОЙЧИВОСТИ  ЛИНЕЙНЫХ ФОРМАЦИЙ }

\author{А.В. Лакеев\footnote{Работа выполнена при  частичной финансовой поддержке Российского фонда фундаментальных исследований (грант \mbox{№\,19-08-00746}).}  \\  
	Институт динамики систем и теории управления   
	им. В.М.\,Матросова \\ СО РАН,  Российская Федерация, 664033, Иркутск, ул.Лермонтова, 134 \\
	e-mail: \texttt{lakeyev@icc.ru}} 

\date{} 

\maketitle 

\begin{abstract}
	Получены необходимые и достаточные условия внутренней устойчивости  формаций, динамика которых 
	определяется линейными дифференциальными уравнениями.
	
	{\it Ключевые слова:} формация, устойчивой по входу-состоянию, ориентированный бесконтурный граф, гурвицева матрица, стабилизируемость.
\end{abstract}

\section{Введение}

Понятие внутренней устойчивости формаций нескольких 
управляемых движущихся объектов возникло при изучении задачи управления группой автономных 
роботов \cite{ba}.

Под формацией понимается конечное множество 
взаимосвязанных движущихся объектов, связи между которыми (типа 
"ведущий-ведомый") задаются некоторым ориентированным бесконтурным
графом. 
Каждой паре связанных объектов сопоставляется некоторый вектор, определяющий 
их требуемое (в идеале) взаимное расположение, которое реально может быть 
выдержано только с некоторой точностью. Изучается задача существования 
управлений, обеспечивающих малые отклонения формации от желаемого идеального 
состояния при малом начальном отклонении и малых возмущающих внешних 
воздействиях. 
Это свойство, названное  внутренней устойчивостью формаций, 
видимо, впервые появилось  в работах H.G. Tanner'а, G.J.~Pappas'а, V.Kumar'а (2002) 
\cite{t1}-\cite{t4}
(см. также \cite{o,l}).

В этих работах за оценку малости отклонения была выбрана концепция "устойчивости от входа к состоянию" (input to state stability). Данное понятие впервые появилось в работах E.D.Sontag’а (1989) \cite{s1,s2}
и активно  исследовалось в работах многих авторов, 
при этом в основном использовался метод функций Ляпунова. Некоторые итоги этих более чем 20--летних исследований содержатся в обзоре \cite{d}.

Отметим, что имеются и различные другие подходы к исследованию поведения формаций. Например, в работах \cite{vku1, vku2, um}, 
также на основе метода функций Ляпунова, изучается свойство типа диссипативности.

Однако практически во всех работах при исследовании внутренней устойчивости формаций вопрос существования управлений, обеспечивающих свойство  внутренней устойчивости формации, в явном виде не ставился.

В данной работе мы будем рассматривать только формации, у которых динамика каждого объекта описывается линейным дифференциальным уравнением и управления ищутся в классе аффинных обратных связей. Работа является расширенным вариантом работы \cite{lak}. 

\section{Основные понятия и определения}

Перейдем к точным формулировкам. Следующие определения взяты из работ \cite{t1,t2,t3,t4} 

	
{\bf О п р е д е л е н и е \   1. (\cite{t1}-\cite{t4})} 	
		
		Будем говорить, что задана формация из $l$ взаимосвязанных 
	объектов, если:
	
	1). Задан ориентированный бесконтурный 
	граф $G=(V,E)$ с множеством вершин 
	$V=\{1,\ldots ,l\}$ и множеством ребер $E\subseteq V\times V$, 
	называемый в дальнейшем графом формаций, и при этом вершины графа  
	отождествляются с движущимися объектами формации, а ребра указывают на 
	наличие взаимосвязей между некоторыми из них.
	Если $(i,j)\in E$, то $i$-ый объект называется ведомым а $j$-ый объект
	--- ведущим (ребро $(i,j)\in E$ считаем ориентированным от $i$ к $j$).
	
	2). Динамика движения $i$-го объекта задается  дифференциальным
	уравнением с управлением
	\begin{equation}\label{(1)}
	{\dot {\rm x}}_i=f_i(t,{\rm x}_i ,\ru_i), \; {\rm x}_i\in R^n, \; \ru_i\in R^m, \; t\in R^1_0=[0,+\infty ).
	\end{equation}
		
	3). Каждому ребру $(i,j)\in E$ сопоставлен вектор $d_{ij}\in R^n$,
	задающий требуемое взаимное расположение (в идеале) $i$-го и $j$-го объектов 
	относительно друг друга, при этом в идеале хотелось бы иметь соотношение 
	${\rm x}_i+d_{ij}={\rm x}_j$.
		
	4). Если динамика движения всех объектов задается линейным  дифференциальным
	уравнением с управлением, то есть 
	${\dot {\rm x}}_i=
	A_i{\rm x}_i + B_i{\rm u}_i
	$, где $A_i\in R^{n\times n}$ и $B_i\in R^{n\times m}$ -- матрицы соответствующих размерностей, то формацию будем называть {\bf линейной}.

{\bf З а м е ч а н и е \   1.} 		
	
	Отметим, что в работах \cite{t1}-\cite{t4}  в качестве графа формации используется обратный граф (\cite[стр. 234]{hr}) $G^{\prime}$ по отношения к графу $G$.  

Будем обозначать $L_i=\{j\in V \, | \; (i,j)\in E\}=\{j_{i1},\ldots ,j_{is_i}\}$~---  
множество ведущих для $i$-го объекта. 
Обозначим также $V_0=\{i\in V \, | \; L_i=\varnothing\}$~--- множество {\bf лидеров} формации, т.е. объекты у которых нет ведущих. 

Если $(i,j)\in E$, 
то обозначим 
${\rm z}_{i,j}={\rm x}_{i}+d_{ij}-{\rm x}_{j}$ -- {\bf отклонение (ошибку)} от идеального состояния между $i$-ым и $j$-ым объектами и ${\rm z}=({\rm z}_{i,j})$ -- $l\times l$-матрица ошибок, 

$
{\rm z}_{i,j}=
\left\{
\begin{array}{lc}
	{\rm x}_{i}+d_{ij}-{\rm x}_{j},&\mbox{если} \; (i,j)\in E,\\
	{\bf 0},&\mbox{если} \; (i,j)\notin E.
\end{array}	
\right. 
$

Для оценки величины отклонений используем стандартные классы функций \cite{h}:
\begin{itemize}
	\item
	${\cal K}$~--- класс функций Хана:  
	$\gamma\in {\cal K},$ если $\gamma:[0,+\infty )\to [0,+\infty )$
	непрерывная, строго возрастающая, $\gamma(0)=0.$
	\item
	функции класса ${\cal K}{\cal L}$:
	$\beta \in {\cal K}{\cal L}$, если $\beta :[0,+\infty )
	\times [0,+\infty )\to [0,+\infty )$, 
	$\beta (\cdot ,t)\in {\cal K}$ при
	фиксированном $t$,  
	$\lim\limits_{t\to \infty}\beta (a,t)=0,$ при
	фиксированном $a\in [0,+\infty )$,
	$\beta (a,t)$ непрерывна по $(a,t)$.
\end{itemize}

Пусть для каждого $i\in V$ задан класс допустимых управлений 
${\cal U}_i$ для $i$-го объекта. 
В общем случае будем считать, что ${\cal U}_i$ -- некоторое подмножество множества функций вида $u_i(t,\rx_i,\rx_{j_{i1}},\ldots ,\rx_{j_{is_i}})$ (то есть $u_i$ является функцией от $t$, $x_i$ и состояний ведущих для $i$-го объекта) таких, что задача Коши для уравнения (\ref{(1)}) с этим управлением имеет единственное решение на $R^1_0$. В дальнейшем мы уточним, какие именно множества функций будут рассматриваться в данной работе.

	{\bf О п р е д е л е н и е \   2.} 	
	
	Формация называется {\bf внутренне устойчивой},  если существуют
	функции $\gamma_i\in {\cal K}, \; i\in V_0 $
	и функция $\beta \in {\cal K}{\cal L}$
	такие, что для любых допустимых управлений
	$u_i\in {\cal U}_i, \; i\in V_0$ найдутся управления $u_j\in {\cal U}_j, \; j\in 
	V\backslash V_0$ такие, что при любых начальных значениях
	${\rm x}_{0i}\in R^n, \; i=\overline{1,l}$ и любых $t\geq 0$ выполняется оценка
	
	\begin{equation}\label{(2)}
	\Vert z(t)\Vert \leq \beta (\Vert z(0)\Vert ,t)+\sum_{i\in V_0}
	\gamma_i(\sup\limits_{\tau \in [0,t]} \Vert u_i \Vert ) ,
	\end{equation}
	где $z(t)$ --- матрица ошибок, полученная на решениях уравнений (\ref{(1)}) при заданных управлениях $u_i\in {\cal U}_i$ и заданных начальных значениях $x_i(0)=\rx_{0i}\in R^n, \; i=\overline{1,l}$.

Здесь и всюду в дальнейшем нормы векторов и матриц считаются евклидовыми.

Для линейных формаций будем рассматривать следующие классы  допустимых управлений:
\begin{itemize}
	\item если $i\in V_0$, то ${\cal U}_i={\cal U}_i^{pr}=
	\{ u_i\; |\;u_i:R_0^1\rightarrow R^m , u_i\;-\;\mbox{непрерывно}  \}$;	
	\item если $i\in V\backslash V_0$ и $L_i=
	\{j_{i1},\ldots ,j_{is_i}\}$, то  
	${\cal U}_i={\cal U}_i^{aff}=\{ u_i \;|\;u_i:R^{n(s_i+1)}\to R^m \} $, 
	
	$u_i({\rm x}_i,{\rm x}_{j_{i1}},\ldots ,{\rm x}_{j_{is_i}})=S_i{\rm x}_i+{\sum_{r=1}^{s_i}}K_{ij_{ir}}{\rm x}_{j_{ir}}+
	k_i=S_i{\rm x}_i+{\sum_{j\in L_i}K_{ij}{\rm x}_{j}}+
	k_i$, 
	
	где  $S_i,K_{ij}\in R^{m\times n}$, $k_i\in R^m$.	
\end{itemize}

{\bf О п р е д е л е н и е  3.} 

Линейную формацию будем называть {\bf линейно внутренне устойчивой  $(\mbox{ЛВУ})$},  если свойство из определения 2 выполняется при ${\cal U}_i={\cal U}_i^{pr}$ для $i\in V_0$, и при ${\cal U}_i={\cal U}_i^{aff}$  для $i\in V\backslash V_0$.

В кванторном виде это свойство более кратко можно записать следующим образом: 

$
\begin{array}{c}
	(\mbox{ЛВУ})\equiv(\exists \gamma_i\in {\cal K}: i\in V_0) (\exists \beta \in {\cal K}{\cal L}) (\forall u_i\in {\cal U}_i^{pr}: i\in V_0)  (\exists u_j\in {\cal U}_j^{aff}: j\in V\backslash V_0) 
	\\
	(\forall {\rm x}_{0i}\in R^n: i\in V) 
	(\forall t\geq 0)\; 
	(\Vert z(t)\Vert \leq \beta (\Vert z(0)\Vert ,t)+\sum_{i\in L_0}
	\gamma_i(\sup\limits_{\tau \in [0,t]} \Vert u_i(\tau) \Vert )).
\end{array}
$

Будем рассматривать еще два свойства линейных формаций:

$(\mbox{ЛВУ})\,_1$ -- получается из $\mbox{(ЛВУ)}$ перестановкой 
кванторов
$(\forall u_i\in {\cal U}_i^{pr}: i\in V_0)$ и 

\noindent
$(\exists u_j\in {\cal U}_j^{aff}: j\in V\backslash V_0)$ или, в кванторном виде,

$
\begin{array}{c}
	\mbox{(ЛВУ)}_1\equiv (\exists \gamma_i\in {\cal K}: i\in V_0) (\exists \beta \in {\cal K}{\cal L})(\exists u_j\in {\cal U}_j^{aff}: j\in V\backslash V_0) 
	(\forall u_i\in {\cal U}_i^{pr}: i\in V_0)
	\\
	(\forall {\rm x}_{0i}\in R^n: i\in V) (\forall t\geq 0)
	\Vert z(t)\Vert \leq \beta (\Vert z(0)\Vert ,t)+
	\sum_{i\in V_0}
	\gamma_i(\sup\limits_{\tau \in [0,t]} \Vert u_i(\tau)\Vert );
\end{array}
$

$(\mbox{ЛВУ})\,_0$ -- получается из $\mbox{(ЛВУ)}$ при нулевых 
управлениях лидеров, т.е. $u_i(t)\equiv {\bf 0}$ для всех $i\in V_0$
(при этом в неравенстве (2) второе слагаемое тождественно равно нулю), 
или, в кванторном виде, 

$
\begin{array}{c}
	\mbox{(ЛВУ)}\,_0\equiv (\exists \beta \in {\cal K}{\cal L})(\exists u_j\in {\cal U}_j^{aff}: j\in V\backslash V_0)  
	(\forall u_i\equiv {\bf 0}: i\in V_0)
	\\
	(\forall {\rm x}_{0i}\in R^n: i\in V) (\forall t\geq 0)
	\Vert z^0(t)\Vert \leq \beta (\Vert z^0(0)\Vert ,t),
\end{array}
$

\noindent
где $z^0(t)$ --- матрица ошибок, полученная на решениях уравнений (\ref{(1)}) при заданных управлениях $u_j\in {\cal U}_j^{aff}$ для $j\in V\backslash V_0$, $u_i(t)\equiv {\bf 0}$ для $i\in V_0$  и заданных начальных значениях $x_i(0)={\rm x}_{0i}\in R^n,$  $ \; i=\overline{1,l}$.

Очевидно, что $(\mbox{ЛВУ})\,_1\Rightarrow (\mbox{ЛВУ})\Rightarrow (\mbox{ЛВУ})\,_0$ .

Рассмотрим более подробно граф формации, при этом будем 
придерживаться терминологии из (\cite[Глава 16]{hr}). 

Во-первых, нетрудно показать, что внутренняя устойчивость формации эквивалентна внутренней устойчивости всех подформаций, определяемых компонентами слабой связности графа $G$. Поэтому в дальнейшем
без ограничения общности будем считать, что граф $G$  -- {\bf слабо связный}. 

Во-вторых, хорошо известно, что ориентированный 
граф $G=(V,E)$ является  бесконтурным тогда и только тогда, когда он допускает перенумерацию вершин с нижней строго треугольной матрицей смежности (\cite[стр. 235, теорема 16.3]{hr}, \cite[стр. 241, теорема 4.8]{nf}).

Более точно, следуя  \cite[стр. 241]{nf}, определим множества уровней ориентированного 
бесконтурного графа $G=(V,E)$, полагая: 

\begin{itemize}
\item $V_0=\{i\in V \, | \; L_i=\varnothing\}$~--- множество лидеров формации;
\item если $V_k$ -- определено и 
$V\setminus\bigcup\limits_{s=0}^k V_s\neq \varnothing$, то 
$V_{k+1}=\{i\in V\setminus\bigcup\limits_{s=0}^k V_s \; | \;L_i\subseteq\bigcup\limits_{s=0}^k V_s \}$;
\item $V_{r}=\{i\in V\setminus\bigcup\limits_{s=0}^{r-1} V_{s} \; | \;L_i\subseteq\bigcup\limits_{s=0}^{r-1} V_s \}$, где $r$ -- наименьшее число, такое, что $V\setminus\bigcup\limits_{s=0}^{r} V_{s}=\varnothing$, то есть
$V=\bigcup\limits_{s=0}^{r} V_{s}$.
\end{itemize}

Множества $V_0, V_1, \ldots , V_r$ называются уровнями графа $G=(V,E)$ и по теореме 4.8 \cite[стр. 241]{nf} для бесконтурного графа являются непустыми множествами, образующими разбиение множества вершин $V$.

Используя это разбиение, 
определим  функцию (так называемую порядковую функцию графа $G$) $O:V\rightarrow \{0,1,\ldots ,r \}$, полагая $O(i)=k$, если $i\in V_k$. Значение  функции  $O(i)$ однозначно определяется по графу $G$.

Если теперь перенумеровать вершины графа  $G$ следующим образом: сначала нумеруем все вершины из $V_0$ (в произвольном порядке), затем все вершины из $V_1$,  затем все вершины из $V_2$  и так далее, то получим нумерацию, {\bf согласованную с порядковой функцией}, то есть для нее выполняется следующее свойство:
\begin{itemize}
\item если $(i,j)\in E$, то $O(i)>O(j)$.	
\end{itemize}
Всюду в дальнейшем будем считать, что изначальная нумерация вершин бесконтурного графа  $G$ {\bf согласована с его порядковой функцией}. Отметим также, что если 
нумерация вершин графа  $G$  согласована с порядковой функцией, то из 
$(i,j)\in E$ следует $i>j$, а множество лидеров $V_0=\{1,\ldots,l_0\}, \;l_0<l$. 

В дальнейшем нам также понадобится функция $p:V\to V$, выбирающая для каждого $i>l_0$ одного из его ведущих. Для определенности будем считать, что это ведущий с наибольшим номером, а если $i\leq l_0$, то  $p(i)=i$, то  есть  

$
p(i)=
\left\{
\begin{array}{lc}
i,&\mbox{если} \; i\leq l_0,\\
\max\{j\;|\;j\in L_i\},&\mbox{если} \; i>l_0.
\end{array}	
\right. 
$

Отметим, что если $i\in V_k,\;k>1$, то найдется $j\in L_i$ такое, что 
$j\in V_{k-1}$, так как иначе $L_i\subseteq \bigcup\limits_{s=0}^{k-2} V_s$ и $O(i)\leq k-1$. Поэтому $p(i)\in V_{k-1}$ и $O(p(i))=O(i)-1$.

Используя итерации функции $p$, получаем последовательность вершин
\linebreak
$i,p(i),p^2(i),\ldots,p^{O(i)}(i)$  
образующих путь  в графе $G$ 
из вершины $i$ в некоторую вершину из $V_0$, которую в дальнейшем будем обозначать  $l_i=p^{O(i)}(i)\in V_0$. 

Также нам в дальнейшем понадобятся следующие вектора $D_i\in R^n$, определяемые рекуррентным образом:

$
D_i=
\left\{
\begin{array}{lc}
{\bf 0},&\mbox{если} \; i\leq l_0,\\
d_{ip(i)}+D_{p(i)},&\mbox{если} \; i>l_0.
\end{array}	
\right. 
$
 
Из рекуррентной формулы для $D_i$ получаем явную формулу  
$D_i=\sum\limits_{s=0}^{O(i)-1}d_{p^{s}(i)p^{s+1}(i)}$.

\section{Необходимые условия для свойства $(\mbox{ЛВУ})\,_0$} 

Путь для формации выполняется свойство $(\mbox{ЛВУ})\,_0$. Выберем и зафиксируем до конца данного раздела функцию $\beta \in {\cal K}{\cal L}$ и  управления для 
объектов с номерами из $V\backslash V_0$, обеспечивающих это свойство, то есть 
матрицы $S_i\in R^{m\times n}$, наборы матриц 
$\{K_{ij}\}_{j\in L_i}, \; K_{ij}\in R^{m\times n}$ и 
вектора $k_i\in R^m$ такие, что если обозначить $x^0_i(t,{\rm x}_{01},\ldots,{\rm x}_{0i})$ -- решения задачи Коши для 
следующей системы дифференциальных уравнений:
\begin{equation}\label{(3)}
	\left\{
	\begin{array}{ll}
		{\dot {\rm x}}_i=A_i{\rm x}_i+B_iu_i(t) , &\mbox{если} \; i\leq l_0,\\
		{\dot {\rm x}}_i=(A_i +B_i S_i){\rm x}_i+B_i({\sum\limits_{j\in L_i}K_{ij}{\rm x}_{j}}+k_i),&\mbox{если} \; i>l_0 ,
	\end{array}	
	\right. 
\end{equation}
с управлениями $u_i(t)\equiv {\bf 0}$ при $i\leq l_0$  и
с начальными условиями 
$x^0_i(0,{\rm x}_{01},\ldots,{\rm x}_{0i})=$
\linebreak
$={\rm x}_{0i},\;i=\overline{1,l}$, и, соответствующие им 
ошибки, 
$z^0_{i,j}(t)= z^0_{i,j}(t,{\rm x}_{01},\ldots,{\rm x}_{0i})=$
\linebreak
$=x^0_i(t,{\rm x}_{01},\ldots,{\rm x}_{0i})-x^0_j(t,{\rm x}_{01},\ldots,{\rm x}_{0j})+d_{ij}$, при $(i,j)\in E$, то для всех 
$t\geq 0$ выполняется оценка
\begin{equation}\label{(4)}
	\Vert z^0(t)\Vert \leq \beta (\Vert z^0(0)\Vert ,t).
\end{equation}

Из неравенства 
(\ref{(4)}) 
получаем следующее необходимое условие.

{\bf Необходимое условие 1.} $\lim\limits_{t\to \infty}z^0_{i,j}(t)=0$ при любых 
$x_i(0)=\rx_{0i}\in R^n, \; i=\overline{1,l}$.

Рассмотрим теперь поведение решений системы (\ref{(3)}) при нулевых начальных условиях. Обозначим эти решения $\tilde{x}^0_i(t)=x^0_i(t,{\bf 0},\ldots,{\bf 0})$. Так как очевидно, что  $\tilde{x}^0_i(t)\equiv{\bf 0}$ при 
$i\leq l_0$,
то индукцией по $O(i)$ из необходимого условия 1 получаем следующее утверждение. 
	
{\bf П р е д л о ж е н и е \   1.} 		
	
	Если для формации выполняется свойство $(\mbox{ЛВУ})\,_0$, то 
\begin{itemize}
	\item[1)]
	для всех $i\in V$ существует конечный предел 
	$\lim\limits_{t\to \infty}\tilde{x}^0_i(t)=-D_i$;
	\item[2)]
	если $j\in L_i$, то
	$\lim\limits_{t\to \infty}\tilde{x}^0_i(t)=-(d_{ij}+D_j)$ и, следовательно, $d_{ij}=D_i-D_j$. 
\end{itemize}
Доказательство этого и последующих утверждений приведено в Приложении. 
	
{\bf З а м е ч а н и е \   2.} 
	
Отметим, что из этого предложения получаем жесткое ограничение на $d_{ij}$ в случае $|L_i|>1$. Если же $|L_i|=1$, то $L_i=\{p(i)\}$ и равенство $d_{ip(i)}=D_i-D_{p(i)}$ выполняется всегда.

Получим теперь еше одно необходимое условие для того, чтобы выполнялось свойство $(\mbox{ЛВУ})\,_0$. Для этого при любом $\rx_{0}\in R^n$ обозначим $\rx_{0i}^*=\rx_{0}-D_i$, $x_i^*(t,\rx_{0})=x^0_i(t,\rx_{01}^*,\ldots,\rx_{0i}^*)$,  
$z^*_{i,j}(t,\rx_{0})=x^*_i(t,\rx_{0})-x^*_j(t,\rx_{0})+d_{ij}$, 
$i=\overline{1,l},\; (i,j)\in E$. Тогда из второго утверждения предложения 1 получаем, что 
$z^*_{i,j}(0,\rx_{0})=\rx_{0}-D_i-(\rx_{0}-D_j)+d_{ij}={\bf 0}$ и поэтому из 
неравенства (\ref{(4)}) получаем следующее тождество. 

{\bf Необходимое условие 2.} $z^*_{i,j}(t,\rx_{0})\equiv {\bf 0}$, при любых 
$\rx_{0}\in R^n, \; t\geq 0, \; i,j=\overline{1,l}$.

Дифференцируя это тождество по $t$ и после этого подставляя $t=0$, получаем следующее утверждение. 

\pagebreak
{\bf П р е д л о ж е н и е \   2.} 	 
	
Если для формации выполняется свойство $(\mbox{ЛВУ})\,_0$, то для любого ребра $(i,j)\in E$
матрицы $S_i,S_j\in R^{m\times n}$, вектора $k_i,k_j\in R^m$ и наборы матриц $\{K_{is}\}_{s\in L_i},\{K_{j\nu}\}_{\nu\in L_j} \; K_{is},K_{j\nu}\in R^{m\times n}$ удовлетворяют следующей системе линейных уравнений: 
$$	
A_i +B_i(S_i+\sum\limits_{s\in L_i}K_{is})=
A_j +B_j(S_i+\sum\limits_{\nu\in L_j}K_{j\nu}),
$$
$$
B_i(k_i-S_iD_i-\sum\limits_{s\in L_i}K_{is}D_s)-A_iD_i=
B_j(k_j-S_jD_j-\sum\limits_{\nu\in L_j}K_{j\nu}D_\nu)-A_jD_j,
$$
где для $j\in V_0$ (в этом случае $L_j=\varnothing$) считаем, что $S_j={\bf 0}$,  
$\sum\limits_{\nu\in L_j}K_{j\nu}={\bf 0}$,  $k_j={\bf 0}$.

Получим теперь некоторые следствия из предложения 2.  Для этого обозначим 
$M_i=A_i +B_i(S_i+\sum\limits_{s\in L_i}K_{is})\in R^{n\times n}$ и
$m_i=B_i(k_i-S_iD_i-\sum\limits_{s\in L_i}K_{is}D_s)-A_iD_i\in R^n$, $i=\overline{1,l}$. Из предложения 2 получаем, что если вершины $i$ и $j$ графа $G$ соединены ребром, то есть $(i,j)\in E$ или $(j,i)\in E$, то $M_i=M_j$ и $m_i=m_j$. Но тогда если существует полупуть (\cite[Глава 16]{hr}) из вершины $i$ в вершину  $j$, то $M_i=M_j$ и $m_i=m_j$. А так как мы считаем, что граф $G$ -- слабо связный, то $M_i=M_j$ и $m_i=m_j$ для любых вершин $i$ и $j$ графа $G$, то есть $M_i$ и $m_i$ не зависят от $i$. 

Заметим, что если $i\in V_0$, то $M_i=A_i$ и $m_i={\bf 0}$, а так как всегда 
$1\in V_0$, то $A_i=A_1$ для  всех  $i\in V_0$.
Из предыдущих рассуждений и предложения 2 
получаем, что  верно следующее утверждение.

{\bf П р е д л о ж е н и е \   3.} 	
	 
Если для формации выполняется свойство $(\mbox{ЛВУ})\,_0$, то
		\begin{itemize}
		\item[1)]
		$A_i=A_1$ для всех $i\leq l_0$, то есть все матрицы  $A_i$ одинаковы при $i\in V_0$;
		\item[2)]
		матрицы $S_i\in R^{m\times n}$, вектора $k_i\in R^m$ и наборы матриц $\{K_{is}\}_{s\in L_i}\; K_{is}\in R^{m\times n}$ удовлтворяют следующей системе линейных уравнений: 
		\begin{equation}\label{(5)}
		B_i(S_i+\sum\limits_{s\in L_i}K_{is})=A_1 -A_i,
		\end{equation}
		\begin{equation}\label{(6)}
		B_i(k_i-S_iD_i-\sum\limits_{s\in L_i}K_{is}D_s)=A_iD_i,
		\end{equation}
		для всех $i>l_0$.
	\end{itemize}

Следующее утверждение касается свойств матриц $\tilde{A}_i=A_i+B_iS_i$ и $A_1$.

	{\bf П р е д л о ж е н и е \   4.} 	
	
	Если для формации выполняется свойство $(\mbox{ЛВУ})\,_0$, то 
\begin{itemize}
		\item[1)]
		для всех $i>l_0$ матрицы $\tilde{A}_i=A_i+B_iS_i$ -- гурвицевы и, следовательно, пары матриц $(A_i,B_i)$ -- стабилизируемы;
		\item[2)]
		если $l_0>1$, то матрица $A_1$ -- гурвицева.
\end{itemize}	

{\bf З а м е ч а н и е \   3.} 	
	
Сделаем несколько замечаний о разрешимости уравнений (\ref{(5)})-(\ref{(6)}). 
	Обозначая $N_i=S_i+{\sum\limits_{s\in L_i}K_{is}}$ и $\tilde{k}_i=k_i-S_iD_i-\sum\limits_{s\in L_i}K_{is}D_s$, получаем, что 
	$N_i\in R^{m\times n}$ и $\tilde{k}_i\in R^m$ являются решениями следующих уравнений: 
	$$
	B_iN_i=A_1-A_i,
	$$
	$$
	B_i\tilde{k}_i=A_iD_i
	$$
	для всех $i>l_0$.
	
	Обратно, если $N_i\in R^{m\times n}$ и $\tilde{k}_i\in R^m$ являются решениями предыдущих уравнений, 
	$S_i\in R^{m\times n}$	-- любая матрица, а набор матриц $\{K_{is}\}_{s\in L_i}\; K_{is}\in R^{m\times n}$  удовлетворяет соотношению $\sum\limits_{s\in L_i}K_{is}=N_i-S_i$, то, полагая 
	$k_i=\tilde{k}_i+S_iD_i+\sum\limits_{s\in L_i}K_{is}D_s$, 
	получаем решение уравнений (\ref{(5)})-(\ref{(6)}). 

\section{Необходимые и достаточные условия \\ линейной внутренней устойчивости}

Набор условий для наличия свойства $(\mbox{ЛВУ})\,_0$ из предложений 1, 3, 4 
оказался не только необходимым, но и достаточным. Более точно, верно  следующуе утверждение.

	{\bf Т е о р е м а \   1.} 	
	
	Для линейной формации выполняется свойство линейной внутренней устойчивости тогда и только тогда, когда  выполнены следующие условия:
	\begin{itemize}
		\item[1)]
	для всех $i> l_0$ пары матриц $(A_i,B_i)$ -- стабилизируемы, то есть 
	существуют $m\times n$-матрицы $S_i$ 	такие, что	матрицы $\tilde{A}_i=A_i+B_iS_i$  -- гурвицевы;
		\item[2)]
	для всех $i> l_0$ существует  $m\times n$-матрица $N_i$ и вектора $\tilde{k}_i\in R^m$,   удовлетворяющие следующим линейным уравнениям:
	\begin{equation}\label{(7)}
	B_iN_i=A_1-A_i,
	\end{equation}
	\begin{equation}\label{(8)}
	B_i\tilde{k}_i=A_iD_i,
	\end{equation}
		\item[3)]
	если  $(i,j)\in E$,  то $d_{ij}=D_i-D_j$;
		\item[4)]
	если $l_0>1$, то $A_i=A_1$ для всех $i\leq l_0$ и матрица $A_1$ -- гурвицева.	
\end{itemize}

\noindent
При этом, если $S_i$, $N_i$ и $\tilde{k}_i$ удовлетворяют условиям 1), 2) теоремы, то  управления ведомых, обеспечивающие внутреннюю устойчивость, можно выбрать в виде 
$$
 u_i(\rx_1,\ldots ,\rx_i)=S_i\rx_i+{\sum_{s\in L_i}K_{is}\rx_{s}}+k_i,\;\;i> l_0,
$$
где набор матриц $\{K_{is}\}_{s\in L_i}\; K_{is}\in R^{m\times n}$  удовлетворяет соотношению $\sum\limits_{s\in L_i}K_{is}=N_i-S_i$, а $k_i=\tilde{k}_i+S_iD_i+\sum\limits_{s\in L_i}K_{is}D_s$, 
независимо от управлений лидеров, т.е. обеспечивающих и свойство 
$(\mbox{ЛВУ})\,_1$.

Сформулируем ряд следствий из этой теоремы.

Так как все необходимые условия получены из свойства $(\mbox{ЛВУ})\,_0$, то верно следующее утверждение. 

		{\bf С л е д с т в и е \   1.} 
	
	Свойства $(\mbox{ЛВУ})_{\,1}$, $(\mbox{ЛВУ})$ и $(\mbox{ЛВУ})\,_0$  эквивалентны.

{\bf З а м е ч а н и е \   4.} 	
	
Заметим, что одно из управлений получается зависящим только от одного $\rx_{s}, \;s\in V_i$ (например, от $\rx_{p(i)}$) и	$\rx_i$  если выбрать $K_{is}=N_i-S_i$ и $K_{is^\prime}={\bf 0}$ для $s^\prime\neq s$, то есть
$$
u_i(\rx_{s},\rx_i)
=S_i(\rx_i+D_i)+(N_i-S_i)(\rx_{s}+D_{s})+\tilde{k}_i,\;\;i>l_0.
$$

Отметим также, что если перебирать {\bf все} $m\times n$-матрицы $S_i$, стабилизирующие пары матриц $(A_i,B_i)$ ($i> l_0$), то есть такие, что	матрица $A_i+B_i S_i$ -- гурвицева,  а также {\bf все}  $N_i$, $\tilde{k}_i$, являющиеся решениями уравнений (\ref{(7)}) , (\ref{(8)}) и {\bf все} наборы матриц $\{K_{is}\}_{s\in L_i}\; K_{is}\in R^{m\times n}$,  удовлетворящие соотношению $\sum\limits_{s\in L_i}K_{is}=N_i-S_i$, то по формуле 
$$
u_i(\rx_1,\ldots ,\rx_i)=S_i(\rx_i+D_i)+{\sum_{s\in L_i}K_{is}(\rx_{s}+D_s)}+\tilde{k}_i
$$
получим {\bf все} управления, при которых выполняется свойство $(\mbox{ЛВУ})_{\,1}$ (естественно, при выполнении условий 3) и 4) теоремы 1.

Далее заметим, что в случае, когда у формации больше одного лидера, то в теореме 1  условие 1) является следствием условий 2) и 4), так как можно взять $S_i=N_i$ и, кроме того, при этом можно взять 
$K_{is}={\bf 0}\in R^{m\times n}$ для всех $i>l_0$, $s\in L_i$. Тогда получаем следующее утверждение.

		{\bf С л е д с т в и е \   2.} 
	
Для линейной формации, имеющей больше одного лидера, выполняется свойство линейной внутренней устойчивости тогда и только тогда, когда  
выполнены следующие условия:
	\begin{itemize}
	\item[1)]
	$A_i=A_1$ для всех $i\leq l_0$ и матрица $A_1$ -- гурвицева;	
	\item[2)]
	выполняются условия условий 2) и 3) теоремы 1.
\end{itemize}
\noindent
При этом если $N_i$ и $\tilde{k}_i$ удовлетворяют условию 2) теоремы 1, то  управления ведомых, обеспечивающие внутреннюю устойчивость, можно выбрать в виде 
$$
u_i(\rx_i)=N_i(\rx_i+D_i)+\tilde{k}_i,\;\;i> l_0,
$$
независимо от состояний 
ведущих и управлений лидеров.	

	{\bf З а м е ч а н и е \   5.} 
	 
Рассмотрим самую простую среди формаций, имеющих больше одного лидера. Она состоит   из трех объектов: 1-ый и 2-ой -- лидеры и ведущие для 3-го. По следствию 2 необходимыми и достаточными условиями ее линейной внутренней устойчивости будут следующие (учитывая, что для этой формации $D_1=D_2={\bf 0}$, $p(3)=2$ и $D_3=d_{32}$): 

1) $A_2=A_1$  и матрица $A_1$ -- гурвицева; 	

2) существуют решения уравнений $B_3N_3=A_1-A_3$ и  $B_3\tilde{k}_3=A_3d_{32}$;

3) $d_{31}=d_{32}$.	

Однако очевидно, что (необходимое) условие $d_{31}=d_{32}$ невыполнимо и бессмысленно для {\bf реальных} объектов. Причем из леммы 1 можно получить, что подобная коллизия будет возникать в {\bf любой} формации, имеющей больше одного лидера. Еще одним странным (противоречащим интуиции) фактом является возможность обеспечить внутреннюю устойчивость (в случае ее наличия) управлением, не зависящим от состояний ведущих (как в следствии 2).

Таким образом, данная концепция внутренней устойчивости не подходит для многолидерных формаций. Причиной этого, видимо, является то, что мы пытаемся выдержать заданное расположение ведомого объекта одновременно относительно {\bf всех} его ведущих (то есть пытаемся угнаться сразу за несколькими зайцами). Возможно, более приемлемой в данном случае может оказаться  концепция внутренней устойчивости из \cite{t4},  в которой используется некоторое агрегированное отклонение $i$-го объекта от идеального состояния вида $e_i=\sum\limits_{j\in L_i}P_{ij}\rz_{i,j}$, где $P_{ij}$ -- некоторые матрицы-проекторы, и требуется выполнение оценки вида (\ref{(2)}) для вектора, составленного из $e_i, \; i>l_0$, а не всей матрицы $\rz$.

Рассмотрим теперь формации с одним лидером. В этом случае $l_0=1$, $V_0=\{1\}$ и условие 4) в теореме 1 отсутствует, а условие 3) теоремы 1 хотя и накладывает дополнительные ограничения на $d_{ij}$, но уже  не является настолько критическим, чтобы быть практически невыполнимым.  Например, для самой простой (в этом случае) формации из трех объектов, граф которой будет бесконтурным треугольником, то есть $G=(V,E)$, $V=\{1,2,3\}$, $E=\{(2,1),(3,1),(3,2)\}$, имеем $D_1={\bf 0}$, $D_2=d_{21}$, $D_3=d_{32}+d_{21}$, условие 3) теоремы 1 сводится к равенству 
$d_{31}=D_3-D_1=d_{32}+d_{21}$, которое означает, что вектора $d_{31}$, $d_{32}$ и $d_{21}$ образуют треугольник в $R^n$ вполне выполнимо и, в некотором смысле, даже естественно. 

Далее рассмотрим случай формаций, для которых условие 3) теоремы 1 отсутствует. Это формации, граф которых является входящим деревом (\cite[стр. 235]{hr}, то есть слабо связным бесконтурным графом, не имеющим полуконтуров, со стоком. Согласно 
\cite[стр. 236, теорема 16.4$^\prime$]{hr} у формации с таким графом будет один лидер и, кроме того, у всех остальных объектов будет единственный ведущий, то есть $L_i=\{p(i)\}$ при $i>1$. В этом случае условие 3) теоремы 1 всегда выполнено и набор матриц  $\{K_{is}\}_{s\in L_i}$     состоит из одной матрицы $K_{ip(i)}$. 
Поэтому получаем следующее утверждение.
	
	{\bf С л е д с т в и е \   3.}
	
Для линейной формации, граф которой является входящим деревом, 
выполняется свойство линейной внутренней устойчивости тогда и только тогда, когда  
выполняются  условия  1), 2) теоремы 1.

При этом, если $S_i$, $N_i$ и $\tilde{k}_i$ удовлетворяют условиям 1), 2) теоремы 1, то  управления ведомых, обеспечивающие внутреннюю устойчивость, можно выбрать в виде 
$$
u_i(\rx_{p(i)},\rx_i)
=S_i(\rx_i+D_i)+(N_i-S_i)(\rx_{p(i)}+D_{p(i)})+\tilde{k}_i,\;\;i>1,
$$
независимо от управления лидера.

\section{Внутренняя устойчивость подформаций}

В работах \cite{t1}-\cite{t4} также затрагивался вопрос о том, как связаны между собой внутренняя устойчивость формация и внутренняя устойчивость ее подформаций.

Имея критерий внутренней устойчивости, 
мы можем ответить и на  этот вопрос. При этом  в данном разделе мы  будем рассматривать в основном только двухэлементные подформации, то есть пары объектов $i,j\in V$ таких, что $(i,j)\in E$.

Из теоремы 1 получаем, что такая  подформация будет линейно внутренне устойчива тогда и только тогда, когда пара $(A_i,B_i)$  стабилизируема
и разрешимы следующие уравнения:
\begin{equation}\label{(9)}
B_iN_{ij}=A_j-A_i,
\end{equation}
\begin{equation}\label{(10)}
B_i\tilde{k}_{ij}=A_id_{ij} .
\end{equation}

Естественно возникают следующие два вопроса:

1) будет ли формация линейно  внутренне устойчива, если все ее двухэлементные подформации линейно  внутренне устойчивы?

2) будут ли все двухэлементные подформации линейно  внутренне устойчивы, если 
формация линейно  внутренне устойчива?

Оказалось, что в общем случае ответы на оба вопроса отрицательные, как показывают следующие два примера.

{\bf П р и м е р \   1.}	
	
Рассмотрим формацию из трех объектов, граф $G$ которой -- бесконтурный треугольник, $G=(V,E)$, $V=\{1,2,3\}$, $E=\{(2,1),(3,1),(3,2)\}$. Пусть $A$ -- любая гурвицева $n\times n$-матрица, $A_1=A_2=A_3=A$. Тогда очевидно, что пары матриц  $(A_2,B_2)$,   $(A_3,B_3)$ -- стабилизируемы при любых $B_2$, $B_3$, и уравнения (\ref{(9)}) имеют решения $N_{ij}={\bf 0}$ для $i=2,3$. Далее будем считать, что управление -- скалярное ($m=1$). Возьмем любой вектор $d\in R^n,\; d\neq {\bf 0}$ и положим $B_2=B_3=Ad$, $d_{21}=d_{31}=d_{32}=d$. Тогда уравнения (\ref{(10)}) имеют решения $\tilde{k}_{ij}=1$ для всех $(i,j)\in E$ и, следовательно, все двухэлементные подформации линейно  внутренне устойчивы. Однако уже упоминавшееся выше необходимое условие внутренне устойчивости всей формации $d_{31}=d_{32}+d_{21}$ не выполняется.

{\bf П р и м е р \   2.}		
	
Рассмотрим формацию из трех объектов, образующих цепь, 
то есть $G=(V,E)$, $V=\{1,2,3\}$, $E=\{(3,2), (2,1)\}$.	Такой граф будет входящим деревом. Поэтому к этой формации применимо следствие 3 и для ее линейной  внутренней устойчивости необходимо и достаточно выполнение следующих условий:

1)  пары матриц  $(A_2,B_2)$,   $(A_3,B_3)$ -- стабилизируемы;

2) существуют матрицы $N_2,\; N_3$  и вектора $\tilde{k}_{2},\; \tilde{k}_{3}$, удовлетворяющие уравнениям

$B_2N_{2}=A_1-A_2$, $B_2\tilde{k}_{2}=A_2D_{2}$,   $B_3N_{3}=A_1-A_3$,    $B_3\tilde{k}_{3}=A_3D_{3}$.
	
\noindent
При этом  $D_1={\bf 0}$, $D_2=d_{21}$, $D_3=d_{32}+d_{21}$.

Возьмем $n=2$, $m=1$ и следующие матрицы и вектора: 
	
$
A_1=\left(
\begin{array}{cc}
1&0\\
0&2
\end{array}	
\right) ,\;
B_1=\left(
\begin{array}{c}
0\\
0
\end{array}	
\right) , 
$
$
A_2=\left(
\begin{array}{rr}
0&-1\\
-1&1
\end{array}	
\right) ,\;
B_2=\left(
\begin{array}{c}
1\\
1
\end{array}	
\right) ,
$

$ 
A_3=\left(
\begin{array}{rr}
0&-1\\
-2&0
\end{array}	
\right) ,\;
B_3=\left(
\begin{array}{c}
1\\
2
\end{array}	
\right) , 
$
$
d_{21}=\left(
\begin{array}{c}
2\\
1
\end{array}	
\right) ,\;
d_{32}=\left(
\begin{array}{c}
2\\
3
\end{array}	
\right) . 
$	

\noindent
Заметим, что пары матриц $(A_2,B_2)$  и	$(A_3,B_3)$ будут вполне управляемы, так как
\linebreak
$
rank\left[ B_2, A_2 B_2    \right]\ = rank\left[ B_3, A_3 B_3    \right] =2 
$
и, следовательно, стабилизируемы.

Далее нетрудно проверить, что при 
$N_2=N_3=\left(
\begin{array}{cc}
1&1
\end{array}	
\right)$ и $\tilde{k}_2=-1,\; \tilde{k}_3=-4$ выполняются уравнения из условия 2), 
поэтому формация линейно  внутренне устойчива.

Однако для двухэлементной подформации, состоящей из 3-го и 2-го объектов, легко показать, что не существует решения уравнения $B_3\tilde{k}_{32}=A_3d_{32}$ (и уравнения $B_3N_{32}=A_2-A_3$) и, следовательно, эта двухэлементная подформация не будет линейно  внутренне устойчива.


{\bf П р и л о ж е н и е. }

	{\bf Доказательство предложения \   1.}
	
	Докажем первое утверждение. Если $O(i)=0$, то есть $i\in L_0$, то утверждение очевидно. Пусть оно верно для всех $j$ таких, что  $O(j)<O(i)$ и $O(i)>0$. Из необходимого условия 1 получаем 
	
	$\lim\limits_{t\to \infty}z^0_{i,p(i)}(t,{\bf 0},\ldots,{\bf 0})=\lim\limits_{t\to \infty}(\tilde{x}^0_i(t)-\tilde{x}^0_{p(i)}(t)+d_{ip(i)})=0,$
	
	\noindent
	поэтому
	
$\lim\limits_{t\to \infty}\tilde{x}^0_i(t)=\lim\limits_{t\to \infty}(\tilde{x}^0_{p(i)}(t)-d_{ip(i)})=-D_{p(i)}-d_{ip(i)}=-D_i.$

Докажем второе утверждение. Если $j\in L_i$, то	опять из необходимого условия 1 и первого утверждения получаем 

	$\lim\limits_{t\to \infty}z^0_{i,j}(t,{\bf 0},\ldots,{\bf 0})=\lim\limits_{t\to \infty}(\tilde{x}^0_i(t)-\tilde{x}^0_{j}(t)+d_{ij})=0,$

\noindent
поэтому

$\lim\limits_{t\to \infty}\tilde{x}^0_i(t)=\lim\limits_{t\to \infty}(\tilde{x}^0_{j}(t)-d_{ij})=-(D_{j}+d_{ij}).$

Но если предел существует, то он единственный. Следовательно, $d_{ij}=D_i-D_j$.
Предложение 1 
доказано. $\blacksquare$

{\bf Доказательство предложения \   2.}	
	
Дифференцируя тождество из необходимого условия 2 по $t$, получаем

$\dot{z}^*_{i,j}(t,\rx_{0})=\dot{x}^*_i(t,\rx_{0})-\dot{x}^*_j(t,\rx_{0})=
(A_i+B_iS_i) x^*_i(t,\rx_{0})+B_i({\sum_{s\in L_i}K_{is}x^*_s(t,\rx_{0})}+k_i)-$
$-(A_j+B_jS_j) x^*_j(t,\rx_{0})-B_j({\sum_{\nu\in L_j}K_{j\nu}x^*_\nu(t,\rx_{0})}+k_j)
\equiv {\bf 0}.$

Подставляя в полученное тождество $t=0$, имеем 

\noindent
$
(A_i+B_iS_i) (\rx_0-D_i)+B_i(\sum_{s\in L_i}K_{is}(\rx_0-D_s)+k_i)-
(A_j+B_jS_j) (\rx_0-D_j) -$
$-B_j(\sum_{\nu\in L_j}K_{j\nu}(\rx_0-D_\nu)+k_j)
\equiv {\bf 0},
$

\noindent
или, после очевидных преобразований,

\noindent
$
(A_i +B_i(S_i+{\sum_{s\in L_i}K_{is}})-
A_j -B_j(S_j+\sum_{\nu\in L_j}K_{j\nu}))\rx_0+
B_ik_i-A_iD_i -$ 
\linebreak
$ -B_i(S_iD_i+\sum_{s\in L_i}K_{is}D_s)-
B_jk_j+A_jD_j +B_j(S_jD_j+\sum_{\nu\in L_j}K_{j\nu}D_\nu)
\equiv {\bf 0}.
$

Так как последнее тождество выполняется при любых $\rx_{0}\in R^n$, то оно эквивалентно равенствам из предложения 2. $\blacksquare$

Для доказательства гурвицевости матрицы $A_1$ при $l_0>1$ нам потребуется следующее утверждение относительно графа $G$.

	{\bf Л е м м а \   1.}
	
	Если $G$ -- слабо связный граф без контуров и  $|V_0|=l_0>1$, то существуют вершины $i,j,s$  такие, что 
	\begin{itemize}
		\item 
		$j,s\in V_0$ -- лидеры, $j\neq s$; 
		\item 
		существует путь  $i,j_1,\ldots, j_{\zeta}=j$ из $i$ в $j$;
		\item
		существует путь  $i,s_1,\ldots, s_{\eta}=s$ из $i$ в $s$.
	\end{itemize}
	
	{\bf Доказательство леммы \   1.}
	
    Разобъем множесва $V_k, \; k=\overline{0,r}$ на подмножества $V_k^s,\; s=\overline{1,l_0}$, полагая:
    \begin{itemize}
    	\item $V_0^s=\{s\},\; s=\overline{1,l_0}$;
    	\item если $V_k^s$  определены, то 
    	$V_{k+1}^s=\{i\in V_{k+1} \; | \;\exists j\in V_k^s \;\;(i,j)\in E \}$.
    \end{itemize}
    
    Покажем индукцией по $k$,
    что для любого $k=\overline{0,r}$ выполняется равенство 
    $V_{k}=\bigcup\limits_{s=1}^{l_0} V_k^s$. 
    Для $k=0$ это очевидно. 
    
    Пусть $V_{k}=\bigcup\limits_{s=1}^{l_0} V_k^s$. Тогда, 
    если $i\in V_{k+1}$, то $(i,p(i))\in E$ и $p(i)\in V_k$. Поэтому, выбирая $s$ такое, что $p(i)\in V_k^s$, получаем $i\in V_{k+1}^s$. Следовательно, $V_{k+1}=\bigcup\limits_{s=1}^{l_0} V_{k+1}^s$.
    
    Заметим, что если $i\in V_k^s$, то существует путь из $i$ в $s$. 
    
    Поэтому, если для некоторого $k>0$ найдутся $j,s\in V_0,\;j\neq s$ такие, что $V_k^j \cap V_k^s\neq\varnothing$, то для этих $j,s$ и любого  $i\in V_k^j \cap V_k^s$ выполняется требуемое свойство. 
    
    Допустим теперь, что для всех  $k=\overline{0,r}$ и любых $j,s\in V_0,\;j\neq s$ имеем  $V_k^j \cap V_k^s=\varnothing$. Тогда очевидно, что множества $V^s=\bigcup\limits_{k=0}^{r} V_{k}^s,\;s=\overline{1,l_0}$ образуют разбиение $V$ на непустые подмножества. Поэтому из слабой связности графа $G$ и $l_0>1$ следует, что для некоторых $j,s\in V_0,\;j\neq s$ найдется ребро $(i,s_1)\in E$, соединяющее некоторую вершину $i\in V^j$ с вершиной $s_1\in V^s$. Но тогда существует путь из $i$ в $j$ (так как $i\in V^j$) и существует путь из $s_1$ в $s$ (так как $s_1\in V^s$), который вместе с ребром $(i,s_1)\in E$ образует путь из $i$ в $s$.
    Лемма 1 
    доказана. $\blacksquare$

\pagebreak
{\bf Доказательство предложения \   3.}	.
	
Это доказательство 
есть в разделе 3.	

{\bf Доказательство предложения \   4.}	
	
Для доказательства гурвицевости матрицы $\tilde{A}_i=A_i+B_iS_i$ 
выпишем дифференциальное уравнение для $z^0_{i,j}(t)=z^0_{i,j}(t,\rx_{01},\ldots,\rx_{0i})$.

Сначала несколько преобразуем уравнение (\ref{(6)}). Используя (\ref{(5)}) и равенство $d_{ij}=D_i-D_j$, получаем

$
B_ik_i=(A_i+B_iS_i)D_i +B_i(\sum\limits_{s\in L_i}K_{is}D_s)=
(A_1-B_i(\sum\limits_{s\in L_i}K_{is}))D_i+B_i(\sum\limits_{s\in L_i}K_{is}D_s)=
$

$
=A_1D_i-B_i(\sum\limits_{s\in L_i}K_{is}(D_i-D_s))=
A_1D_i-B_i(\sum\limits_{s\in L_i}K_{is}d_{is}),
$
то есть
\begin{equation}\label{(a1)}
B_ik_i=
A_1D_i-B_i(\sum\limits_{s\in L_i}K_{is}d_{is}).
\end{equation}
Далее, дифференцируя $z^0_{i,j}(t)$ и используя (\ref{(5)}), (\ref{(6)}) и (\ref{(a1)}), получаем

\noindent
$
\dot{z}^0_{i,j}(t)=\dot{x}^0_i(t)-\dot{x}^0_j(t)= 
\tilde{A}_ix^0_i(t) +B_i(\sum\limits_{s\in L_i}K_{is}x^0_s(t)+k_i)-
\tilde{A}_jx^0_j(t) -B_i(\sum\limits_{\nu\in L_j}K_{j\nu}x^0_{\nu}(t)+k_j)=
$

\noindent
$
=(A_1-B_i({\sum\limits_{s\in L_i}K_{is}}))x^0_i(t)+B_i(\sum\limits_{s\in L_i}K_{is}x^0_s(t))+A_1D_i-B_i(\sum\limits_{s\in L_i}K_{is}d_{is})-
$

\noindent
$
-(A_1-B_j({\sum\limits_{\nu\in L_j}K_{j\nu}}))x^0_j(t)-
B_j(\sum\limits_{\nu\in L_j}K_{j\nu}x^0_{\nu}(t))-
A_1D_j+B_j(\sum\limits_{\nu\in L_j}K_{j\nu}d_{j\nu})=
$

\noindent
$
=A_1(x^0_i(t)-x^0_j(t)+D_i-D_j)-B_i({\sum\limits_{s\in L_i}K_{is}(x^0_i(t)-x^0_s(t)+d_{is})})+
$

\noindent
$
+B_j({\sum\limits_{\nu\in L_j}K_{j\nu}(x^0_j(t)-x^0_{\nu}(t)+d_{j\nu})})=
A_1z^0_{i,j}(t)-B_i({\sum\limits_{s\in L_i}K_{is}z^0_{i,s}(t)})+
B_j({\sum\limits_{\nu\in L_j}K_{j\nu}z^0_{j,\nu}(t)})
$, или
\begin{equation}\label{(a2)}
\dot{z}^0_{i,j}(t)=
A_1z^0_{i,j}(t)-B_i({\sum\limits_{s\in L_i}K_{is}z^0_{i,s}(t)})+
B_j({\sum\limits_{\nu\in L_j}K_{j\nu}z^0_{j,\nu}(t)}).
\end{equation}

В дальнейшем нам понадобятся некоторые алгебраические соотношения между $z_{i,j}$. Если $(i,j)\in E$, то рассмотрим путь $i,j,p(j),p_2(j),\ldots,p^{O(j)}(j)=l_j$ из вершины $i$ в 
вершину $l_j\in V_0$ и сумму  ошибок вдоль этого пути. При этом будем использовать следующую формулу для ошибки:  $z_{i,j}=\rx_i-\rx_j+d_{ij}=\rx_i-\rx_j+D_i-D_j=\rx_i+D_i-(\rx_j+D_j)$.  
Тогда получаем следующее очевидное равенство:

\noindent
$
z_{i,j}+z_{j,p(j)}+z_{p(j),p^2(j)}+\ldots+z_{p^{O(j)-1}(j),p^{O(j)}(j)}=
\rx_i+D_i-(\rx_j+D_j)+(\rx_j+D_j)-(\rx_{p(j)}+D_{p(j)})+\ldots+(\rx_{p^{O(j)-1}(j)}+D_{p^{O(j)-1}(j)})-(\rx_{p^{O(j)}(j)}+D_{p^{O(j)}(j)})=\rx_i+D_i-(\rx_{p^{O(j)}(j)}+D_{p^{O(j)}(j)})=\rx_i-\rx_{l_j}+D_i.
$

\noindent
Аналогично, если $(i,s)\in E$, то, рассмотривая путь $i,s,p(s),p_2(s),\ldots,p^{O(s)}(s)=l_s$, получаем равенство

\noindent
$
z_{i,s}+z_{s,p(s)}+z_{p(s),p^2(s)}+\ldots+z_{p^{O(s)-1}(s),p^{O(s)}(s)}=
\rx_i-\rx_{l_s}+D_i
$

\noindent
и, следовательно

\noindent
$
z_{i,s}+z_{s,p(s)}+z_{p(s),p^2(s)}+\ldots+z_{p^{O(s)-1}(s),p^{O(s)}(s)}=
\rx_i-\rx_{l_s}+D_i=\rx_i-\rx_{l_j}+D_i+\rx_{l_j}-\rx_{l_s}=z_{i,j}+z_{j,p(j)}+z_{p(j),p^2(j)}+\ldots+z_{p^{O(j)-1}(j),p^{O(j)}(j)}+\rx_{l_j}-\rx_{l_s}
$

\noindent
или 
\begin{equation}\label{(a3)}
z_{i,s}=z_{i,j}+R_{j,s}(z)+\rx_{l_j}-\rx_{l_s},
\end{equation}
где
\begin{equation}\label{(a4)}
\begin{array}{rcl}
R_{j,s}(z)&=&z_{j,p(j)}+z_{p(j),p^2(j)}+\ldots+z_{p^{O(j)-1}(j),p^{O(j)}(j)}-
\\
& &-(z_{s,p(s)}+z_{p(s),p^2(s)}+\ldots+z_{p^{O(s)-1}(s),p^{O(s)}(s)}),
\end{array}
\end{equation}
то есть $R_{j,s}(z)$ -- линейная комбинация некоторых $z_{\zeta,\eta}$ таких, что $O(\eta)<O(\zeta)<O(i)$.

Подставляя (\ref{(a3)}) в (\ref{(a2)}) и используя (\ref{(5)}), получаем

\noindent
$
\dot{z}^0_{i,j}(t)=
A_1z^0_{i,j}(t)-B_i({\sum\limits_{s\in L_i}K_{is}(z^0_{i,j}(t)+R_{j,s}(z^0(t))+x_{l_j}^0(t)-x_{l_s}^0(t))})+
+B_j({\sum\limits_{\nu\in L_j}K_{j\nu}z^0_{j,\nu}(t)})=
(A_1-B_i(\sum\limits_{s\in L_i}K_{is}))z^0_{i,j}(t)-B_i({\sum\limits_{s\in L_i}K_{is}(R_{j,s}(z^0(t))+x_{l_j}^0(t)-x_{l_s}^0(t))})+
+B_j({\sum\limits_{\nu\in L_j}K_{j\nu}z^0_{j,\nu}(t)})=
\tilde{A}_iz^0_{i,j}(t)-B_i({\sum\limits_{s\in L_i}K_{is}(R_{j,s}(z^0(t))+x_{l_j}^0(t)-x_{l_s}^0(t))})+
+B_j({\sum\limits_{\nu\in L_j}K_{j\nu}z^0_{j,\nu}(t)})
$, то есть
\begin{equation}\label{(a5)}
\dot{z}^0_{i,j}(t)=
\tilde{A}_iz^0_{i,j}(t)-B_i({\sum\limits_{s\in L_i}K_{is}(R_{j,s}(z^0(t))+x_{l_j}^0(t)-x_{l_s}^0(t))})+
B_j({\sum\limits_{\nu\in L_j}K_{j\nu}z^0_{j,\nu}(t)}).
\end{equation}

Выберем теперь следующие начальные значения: $\rx^*_{\zeta 0}=\rx_0-D_{\zeta}$ при  $O(\zeta)<O(i)$, $\rx_0,\,\rx_{0 \zeta}\in R^n$ -- любые, при $O(\zeta) \geq  O(i)$.
Тогда $x_{l_j}^0(t,\rx_{l_j}^*)=x_{l_s}^0(t,\rx_{l_s}^*)$, так как $A_{l_j}=A_{l_s}$, и 
$z_{\zeta,\eta}(t,\rx^*_{0 1},\ldots ,\rx^*_{0\zeta})=z^*_{\zeta,\eta}(t,\rx_{0})\equiv {\bf 0}$ при  $O(\zeta)<O(i)$, в силу необходимого условия 2. Поэтому уравнение (\ref{(a5)}) для $z^0_{i,j}(t)=z^0_{i,j}(t,\rx^*_{01},\ldots ,\rx_{0i})$ при этих начальных значениях принимает следующий вид: 
$$
\dot{z}^0_{i,j}(t,\rx^*_{01},\ldots ,\rx_{0i})=
\tilde{A}_iz^0_{i,j}(t,\rx^*_{01},\ldots ,\rx_{0i}),
$$
и поэтому $z^0_{i,j}(t,\rx^*_{01},\ldots ,\rx_{0i})=e^{t\tilde{A}_i}z^0_{i,j}(0,\rx^*_{01},\ldots ,\rx_{0i})=e^{t\tilde{A}_i}(\rx_{0i}-\rx^*_{0j}+d_{ij})=e^{t\tilde{A}_i}(\rx_{0i}-\rx_{0}+D_{i})$. Тогда из необходимого условия 1 получаем, что 
$\lim\limits_{t\to \infty}e^{t\tilde{A}_i}(\rx_{0i}-\rx_{0}+D_{i})=$
\linebreak
$=\lim\limits_{t\to \infty}z^0_{i,j}(t,\rx^*_{01},\ldots ,\rx_{0i})=
{\bf 0}$, 
при любых $\rx_0,\,\rx_{0i}\in R^n$.  Следовательно,  матрица $\tilde{A}_i=A_i+B_iS_i$ -- гурвицева. 

Покажем теперь, что если $|V_0|=l_0>1$, то и $A_1$ -- гурвицева. Согласно лемме 1 
выбираем вершины $i,j,s$ такие, что $j,s\in V_0$, $j\neq s$ и существуют пути  $i,j_1,\ldots, j_{\zeta}=j$ из $i$ в $j$, и $i,s_1,\ldots, s_{\eta}=s$ из $i$ в $s$. Опять рассматривая суммы ошибок вдоль этих путей, получаем (учитывая что 
$D_{j_{\zeta}}=D_j={\bf 0}$, $D_{s_{\eta}}=D_s={\bf 0}$ и $z_{i,j}=\rx_i+D_i-(\rx_{j}+D_{j})$),

\noindent   
$
z_{i,j_1}+z_{j_1,j_2}+\ldots+z_{j_{\zeta-2},j_{\zeta-1}}+z_{j_{\zeta-1},j_{\zeta}}=\rx_i+D_i-(\rx_{j_1}+D_{j_1})+(\rx_{j_1}+D_{j_1})-(\rx_{j_2}+D_{j_2})+\ldots+
(\rx_{j_{\zeta-2}}+D_{j_{\zeta-2}})-(\rx_{j_{\zeta-1}}+D_{j_{\zeta-1}})+
(\rx_{j_{\zeta-1}}+D_{j_{\zeta-1}})-(\rx_{j_{\zeta}}+D_{j_{\zeta}})=\rx_i-\rx_{j}+D_i
$

\noindent 
и, аналогично,

\noindent   
$
z_{i,s_1}+z_{s_1,s_2}+\ldots+z_{s_{\zeta-2},s_{\eta-1}}+z_{s_{\eta-1},s_{\eta}}=\rx_i+D_i-(\rx_{s_1}+D_{s_1})+(\rx_{s_1}+D_{s_1})-(\rx_{s_2}+D_{s_2})+\ldots+
(\rx_{s_{\eta-2}}+D_{s_{\eta-2}})-(\rx_{s_{\eta-1}}+D_{s_{\eta-1}})+
(\rx_{s_{\eta-1}}+D_{s_{\eta-1}})-(\rx_{s_{\eta}}+D_{s_{\eta}})=\rx_i-\rx_{s}+D_i.
$

Из последних двух равенств получаем

\noindent   
$
z_{i,j_1}+z_{j_1,j_2}+\ldots+z_{j_{\zeta-2},j_{\zeta-1}}+z_{j_{\zeta-1},j_{\zeta}}-(z_{i,s_1}+z_{s_1,s_2}+\ldots+z_{s_{\zeta-2},s_{\eta-1}}+z_{s_{\eta-1},s_{\eta}})=
\;\;\;\;\;\;\;\;\;\;\;\;\;\;\;\;\;\;\;\;\;\;\;\;\;\;\;=\rx_i-\rx_{j_{\zeta}}+D_i-(\rx_i-\rx_{s_{\eta}}+D_i)=\rx_{s}-\rx_{j},
$ и, следовательно, 

\noindent   
$
z^0_{i,j_1}(t)+z^0_{j_1,j_2}(t)+\ldots+z^0_{j_{\zeta-2},j_{\zeta-1}}(t)+z^0_{j_{\zeta-1},j_{\zeta}}(t)-(z^0_{i,s_1}(t)+z^0_{s_1,s_2}(t)+\ldots+z^0_{s_{\zeta-2},s_{\eta-1}}(t)+
+z^0_{s_{\eta-1},s_{\eta}}(t))=
x^0_{s}(t)-x^0_{j}(t).
$

Далее, так как $x^0_{j}(t)=e^{tA_1}\rx_{0j}$ и $x^0_{s}(t)=e^{tA_1}\rx_{0s}$, то из этого равенства и необходимого условия 1 следует, что 
$\lim\limits_{t\to \infty}(x^0_{s}(t)-x^0_{j}(t))=
\lim\limits_{t\to \infty}e^{tA_1}(\rx_{0s}-\rx_{0j})=
{\bf 0}$ при любых 
$\rx_{0s},\rx_{0j} \in R^n$. Следовательно,  матрица $A_1$ -- гурвицева. 
Предложение 4 
доказано. $\blacksquare$

{\bf Доказательство теоремы \   1.}	
	
	Необходимость условий теоремы 1 
	следует непосредственно из предложений 1, 3, 4 
	и замечания 3. 
	
	Докажем, что эти условия будут и достаточными. Для этого из условий 1), 2) (для всех $i> l_0$) теоремы 1 
	находим   $m\times n$-матрицы $S_i$  	такие, что	матрицы $\tilde{A}_i=A_i+B_i S_i$ -- гурвицевы, 
	$m\times n$-матрицы $N_i$ и вектора $\tilde{k}_i\in R^m$,   удовлетворяющие  линейным уравнениям (\ref{(7)}), (\ref{(8)}), и 
	набор матриц $\{K_{is}\}_{s\in L_i}\; K_{is}\in R^{m\times n}$,  удовлетворяющий соотношению $\sum\limits_{s\in L_i}K_{is}=N_i-S_i$. Согласно замечанию 3, 
	при этом будут выполняться соотношения (\ref{(5)}), (\ref{(6)}).  
	Далее полагаем 
	$$
	u_i(\rx_1,\ldots ,\rx_i)=S_i\rx_i+{\sum_{s\in L_i}K_{is}\rx_{s}}+k_i,
	$$
	при $i> l_0$, где $k_i=\tilde{k}_i+S_iD_i+\sum\limits_{s\in L_i}K_{is}D_s$, 
	а при $i\leq l_0$ выбираем 
	$ u_i:R_0^1\rightarrow R^m$, $u_i\in {\cal U}_i^{pr}$. Подставляя эти управления в систему (\ref{(3)}) и задавая начальные условия $\rx_{0i},\;i\in V$ при $t=0$, находим 
	решение задачи Коши для системы (\ref{(3)}) $x_i(t)=x_i(t,\rx_{01},\ldots,\rx_{0i},u_1,\ldots,u_i)$ такие, что $x_i(0)=x_i(0,\rx_{01},\ldots,\rx_{0i},u_1,\ldots,u_i)=\rx_{0i}$. После этого для всех $(i,j)\in E$ находим ошибки $z_{i,j}(t)=x_i(t)-x_j(t)+d_{ij}$.
	
	Теперь для доказательствы теоремы достаточно показать, что для всех $(i,j)\in E$ найдутся $C_{i,j},\;\alpha_{i,j},\;\beta_{i,j}>0$ 
	такие, что выполняется неравенство
	\begin{equation}\label{(a6)}
	\Vert z_{i,j}(t)\Vert \leq C_{i,j}e^{-\alpha_{i,j}t} \Vert z(0)\Vert +\beta_{i,j}(\sum_{s=1}^{l_0}
	\sup\limits_{\tau \in [0,t]} \Vert u_s(\tau) \Vert ) .
	\end{equation}
	
	Выпишем сначала дифференциальное уравнение для $z_{i,j}(t)$. Заметим, что 
	в данном случае дифференциальное уравнение для $x_i(t)$ можно записать в виде,
	аналогичном (\ref{(3)}):
	\begin{equation}\label{(a7)}
	{\dot {\rm x}}_i=\tilde{A}_i \rx_i+B_i({\sum_{s\in L_i}K_{is}\rx_{s}}+k_i+\delta_i u_i),\;
	i=\overline{1,l}, 
	\end{equation}
	где $\tilde{A}_i=A_i+B_iS_i$ и 
	будем считать, что при $i\leq l_0$ 
	$S_i={\bf 0}$ -- нулевые $m\times n$-матрицы, $k_i\in R^m, k_i={\bf 0}$ -- нулевые вектора,  сумма   $ {\sum_{s\in L_i}K_{is}\rx_{s}}$ в случае 
	$L_i=\varnothing$ равна нулю, а $\delta_i$ определяется следующим образом:
	
	$
	\delta_i=
	\left\{
	\begin{array}{lc}
	1,&\mbox{если} \; O(i)=0,\\
	0,&\mbox{если} \; O(i)>0.
	\end{array}	
	\right. 
	$
		
	Кроме того, так как по условию 3) теоремы 1 
	$d_{ij}=D_i-D_j$, то выполняется равенство (\ref{(a1)}). Поэтому, дифференцируя $z_{i,j}(t)$ по $t$, получаем (аналогично (\ref{(a2)}))
	\begin{equation}\label{(a8)}
	\dot{z}_{i,j}(t)=
	A_1z_{i,j}(t)-B_i({\sum\limits_{s\in L_i}K_{is}z_{i,s}(t)})+
	B_j({\sum\limits_{\nu\in L_j}K_{j\nu}z_{j,\nu}(t)})-\delta_jB_j u_j(t).
	\end{equation}
	
	Далее, опять учитывая, что $d_{ij}=D_i-D_j$, получаем, что для $z_{i,s}(t)$ выполняются алгебраические соотношения (\ref{(a3)}). Подставляя эти соотношения в (\ref{(a8)}) (аналогично (\ref{(a5)})), получаем следующее дифференциальное уравнение:  
	$$
	\begin{array}{rcl}
	\dot{z}_{i,j}(t) & = &
	\tilde{A}_iz_{i,j}(t)-B_i(\sum\limits_{s\in L_i}K_{is}R_{j,s}(z(t)))-
	B_i(\sum\limits_{s\in L_i}K_{is}(x_{l_j}(t)-x_{l_s}(t)))+
	\\
	&  &+B_j({\sum\limits_{\nu\in L_j}K_{j\nu}z_{j,\nu}(t)})-\delta_jB_j u_j(t).
	\end{array}
	$$
    Но, в отличие от (\ref{(a5)}), в нем есть слагаемые вида $x_{l_j}(t)-x_{l_s}(t)$, которые, используя уравнение (\ref{(a7)}) и формулу Коши, можно преобразовать следующим образом:
    
    $
    x_{l_j}(t)-x_{l_s}(t)=e^{tA_1}(\rx_{0l_j}-\rx_{0l_s})+
    \int\limits^{t}_{0}e^{(t-\tau)A_1}(u_{l_j}(\tau)-u_{l_s}(\tau)) d\tau.
    $

    Кроме того, опять из алгебраических соотношений (\ref{(a3)}), но при  $t=0$, получаем 

$
x_{l_j}(t)-x_{l_s}(t)=e^{tA_1}(z_{i,s}(0)-z_{i,j}(0)-R_{j,s}(z(0)))+
\int\limits^{t}_{0}e^{(t-\tau)A_1}(u_{l_j}(\tau)-u_{l_s}(\tau)) d\tau.
$

Из этого соотношения и линейности $R_{j,s}(z)$ получаем следующее дифференциальное уравнение для $z_{i,j}(t)$
	\begin{equation}\label{(a9)}
    \begin{array}{l}
    \dot{z}_{i,j}(t)  = 
    \tilde{A}_iz_{i,j}(t)-B_i(\sum\limits_{s\in L_i}K_{is}R_{j,s}(z(t)-e^{tA_1} \circledast z(0)))
    +B_j({\sum\limits_{\nu\in L_j}K_{j\nu}z_{j,\nu}(t)})+
    \\
   \!+\! B_i(\sum\limits_{s\in L_i}K_{is}(e^{tA_1}(z_{i,j}(0)-z_{i,s}(0))+\int\limits^{t}_{0}e^{(t-\tau)A_1}(u_{l_s}(\tau)-u_{l_j}(\tau)) d\tau))-\delta_jB_j u_j(t),
    \end{array}
	\end{equation}
где для $l\times l$-матрицы $\rz=(\rz_{i,j})$ произведение 
$e^{tA_1} \circledast \rz$ определяется следующим образом:  
$e^{tA_1} \circledast \rz = (e^{tA_1} \rz_{i,j})$. 

Тогда по формуле Коши получаем
\begin{equation}\label{(a10)}
\begin{array}{c}
z_{i,j}(t)  = 
e^{t\tilde{A}_i}z_{i,j}(0)-\int\limits^{t}_{0}e^{(t-\tau)\tilde{A}_i}B_i(\sum\limits_{s\in L_i}K_{is}R_{j,s}(z(\tau)-e^{\tau A_1}\circledast z(0)))d\tau+
\\
+\int\limits^{t}_{0}e^{(t-\tau)\tilde{A}_i}B_j({\sum\limits_{\nu\in L_j}K_{j\nu}z_{j,\nu}(\tau)})d\tau-\int\limits^{t}_{0}e^{(t-\tau)\tilde{A}_i}\delta_jB_j u_j(\tau)d\tau+
\\
\!+\!\int\limits^{t}_{0}e^{(t-\tau)\tilde{A}_i} B_i(\sum\limits_{s\in L_i}K_{is}(e^{\tau A_1}(z_{i,j}(0)-z_{i,s}(0))\!+\!\int\limits^{\tau}_{0}e^{(\tau-\mu)A_1}(u_{l_s}(\mu)-u_{l_j}(\mu)) d\mu))d\tau.
\end{array}
\end{equation}

Покажем теперь, что выполняется оценка (\ref{(a6)}) индукцией по $O(i)$. При этом мы используем хорошо известную оценку матричной экспоненты 
(\cite[стр. 43]{dk}), из которой, в частности, следует, что если матрица $A$ -- гурвицева, то {\bf для любого $\alpha >0$} такого, что $\max \{ {\rm Re}\ \lambda\;|\;\lambda\in \sigma(A) \}<-\alpha$, существует $C=C(\alpha)>0$ такое, что $\Vert e^{tA}\Vert\leq Ce^{-\alpha t}$ при всех $t\geq 0$ (где $\sigma(A)$ -- спектр матрицы $A$).

Пусть $O(i)=1$. Тогда $O(j)=0$ и в формуле (\ref{(a9)}) (и, соответственно, в (\ref{(a10)})) $L_j=\varnothing$, $\delta_j=1$ и  $R_{j,s}(z)={\bf 0}$, для всех $s\in L_i$. Кроме того, если $|L_i|=1$, то последнее слагаемое в  (\ref{(a10)}) равно нулю (так как в сумме $\sum\limits_{s\in L_i}$ присутствует только одно слагаемое с $s=j$), а если $|L_i|>1$, то матрица $A_1$ -- гурвицева. Тогда формула (\ref{(a10)}) принимает следующий вид: 
\begin{equation}\label{(a11)}
\begin{array}{c}
z_{i,j}(t)  = 
e^{t\tilde{A}_i}z_{i,j}(0)-\int\limits^{t}_{0}e^{(t-\tau)\tilde{A}_i}B_j u_j(\tau)d\tau+
\\
\!+\!\int\limits^{t}_{0}e^{(t-\tau)\tilde{A}_i} B_i(\sum\limits_{s\in L_i}K_{is}(e^{\tau A_1}(z_{i,j}(0)-z_{i,s}(0))\!+\!\int\limits^{\tau}_{0}e^{(\tau-\mu)A_1}(u_{l_s}(\mu)-u_{l_j}(\mu)) d\mu))d\tau.
\end{array}
\end{equation}

Так как по условиям 1) и 4) теоремы 1 матрицы $\tilde{A}_i$ и $A_1$ -- гурвицевы (для $A_1$ -- если $|L_i|>1$), то для них выполняются следующие оценки: 

$\Vert e^{t \tilde{A}_i}\Vert\leq \tilde{C}_ie^{-\tilde{\alpha}_i t}$, $\Vert e^{tA_1}\Vert\leq C_1e^{-\alpha_1 t}$  

\noindent
при всех $t\geq 0$ и, кроме того, можно считать, что $\tilde{\alpha}_i<\alpha_1$.

Тогда оценка (\ref{(a6)}) для $z_{i,j}(t)$ при $O(i)=1$  получается в виде соответсвующей суммы из следующих легко проверяемых оценок:  

$
\Vert e^{t \tilde{A}_i}z_{i,j}(0)\Vert\leq \tilde{C}_ie^{-\tilde{\alpha}_i t}\Vert z_{i,j}(0)\Vert ,
$

$
\Vert
\int\limits^{t}_{0}e^{(t-\tau)\tilde{A}_i}B_j u_j(\tau)d\tau
\Vert \leq 
\tilde{\alpha}_i^{-1}\tilde{C}_i\Vert B_j\Vert \sup\limits_{\tau \in [0,t]} \Vert u_j(\tau) \Vert ,
$

$
\Vert
\int\limits^{t}_{0}e^{(t-\tau)\tilde{A}_i} B_iK_{is}e^{\tau A_1}z_{i,s}(0)d\tau
\Vert\leq
(\alpha_1-\tilde{\alpha}_i)^{-1}\tilde{C}_iC_1\Vert B_i\Vert \Vert K_{is}\Vert e^{-\tilde{\alpha}_i t}\Vert z_{i,s}(0)\Vert ,
$

$
\Vert
\int\limits^{t}_{0}e^{(t-\tau)\tilde{A}_i} B_iK_{is}
(\int\limits^{\tau}_{0}e^{(\tau-\mu)A_1}u_{l_s}(\mu)d\mu ) d\tau
\Vert\leq
(\alpha_1\tilde{\alpha}_i)^{-1}\tilde{C}_iC_1\Vert B_i\Vert \Vert K_{is}\Vert \sup\limits_{\tau \in [0,t]} \Vert u_{l_s}(\tau) \Vert .
$

Допустим теперь, что для всех $z_{\zeta,\eta}$ таких, что $O(\eta)<O(\zeta)<O(i)$, оценка типа (\ref{(a6)}) -- выполняется и покажем, что она  выполняется и для $z_{i,j}(t)$. 

Для этого достаточно показать, что каждое слагаемое в формуле (\ref{(a10)}) удовлетворяет соответствующей оценке. Но для всех слагаемых, кроме второго, это уже показано выше, а для второго, учитывая, что  
$R_{j,s}(z)$ -- линейная комбинация некоторых $z_{\zeta,\eta}$ таких, что $O(\eta)<O(\zeta)<O(i)$, достаточно заметить, что если 

$
\Vert z_{\zeta,\eta}(t)\Vert \leq C_{\zeta,\eta}e^{-\alpha_{\zeta,\eta}t} \Vert z(0)\Vert +\beta_{\zeta,\eta}(\sum_{s=1}^{l_0}
\sup\limits_{\tau \in [0,t]} \Vert u_s(\tau) \Vert ) 
$

\noindent
и $\tilde{\alpha}_i<\alpha_{\zeta,\eta}$ (что всегда можно считать выполненным),
то для любой матрицы $B$ выполняется легко проверяемое неравенство 

$
\Vert
\int\limits^{t}_{0}e^{(t-\tau)\tilde{A}_i}Bz_{\zeta,\eta}(\tau) d\tau
\Vert\leq
(\alpha_{\zeta,\eta}-\tilde{\alpha}_i)^{-1}\tilde{C}_iC_{\zeta,\eta}
\Vert B\Vert e^{-\tilde{\alpha}_i t}  \Vert z(0)\Vert + 
$

$
\;\;\;\;\;\;\;\;\;\;\;\;\;\;\;\;\;\;\;\;\;\;\;\;\;\;\;\;\;\;\;\;\;\;\;\;\;\;\;\;\;\;\;\;\;\;\;\;\;\;\;\;\;\;
+\tilde{\alpha}_i^{-1}\tilde{C}_i
\Vert B\Vert
\beta_{\zeta,\eta}(\sum_{s=1}^{l_0}
\sup\limits_{\tau \in [0,t]} \Vert u_s(\tau) \Vert ).
$

\noindent
Опять суммируя соответствующие оценки, получаем, что и для $z_{i,j}(t)$ 
выполняется оценка (\ref{(a6)}).
Теорема 1 
доказана. $\blacksquare$

\section{Заключение}

В работе для  формаций движущихся объектов, динамика которых описывается линейными дифференциальными уравнениями с управлением, получены необходимые и достаточные условия внутренней устойчивости. При этом в качестве классов допустимых управлений для лидеров выбраны программные управления, а для  объектов, имеющих ведущих -- аффинные обратные связи, зависящие от состояния самого объекта и сотояний его ведущих. Полученные условия легко проверяемы и состоят из требований стабилизируемости пары матриц для уравнений ведомых объектов, гурвицевости и совпадения матриц для лидеров (в случае многолидерности), разрешимости некоторых линейных уравнений и ограничений типа равенств на вектора, задающие требуемое взаимное расположение 
между ведомым и ведущим. Кроме  того, описан весь класс управлений, обеспечивающих выполнение свойства линейной  внутренней устойчивости.

Отмечено, что  условия типа равенств   показывают практическую невыполнимость и бессмысленность рассматриваемого понятия внутренней устойчивости  для  
{\bf реальных} 
\linebreak
объектов в случае {\bf нескольких лидеров}. Для формаций с одним лидером выделен класс (формации, граф которых  является входящим деревом), для которого не возникает ограничений типа равенств.  

Исследован вопрос о возможности 
проверки наличия внутренней устойчивости формации, используя ее двухэлементные подформации. Однако на примерах показано, что в общем случае существуют формации (даже с одним лидером), не являющиеся внутренне устойчивыми, у которых все двухэлементные подформации -- внутренне устойчивы, и наоборот, существуют внутренне устойчивые формации (граф которых  является входящим деревом), у которых не все двухэлементные подформации -- внутренне устойчивы.


Lakeyev A. V.    Criterion of Internal Stability of Linear Formations

Matrosov Institute for System Dynamics and Control Theory, Siberian Branch of Russian Academy of Sciences,
Irkutsk, Russian Federation

The necessary and sufficient conditions of internal stability for the formations, whose dy\-na\-mics is defined by linear differential equations, have been obtained. In this case, programmed controls have been chosen in the capacity of the classes of admissible controls for the leaders; while in the case of objects having leaders for them – affine feedbacks dependent on the state of the object itself and on the states of its leaders have been chosen. The conditions obtained may be easily verified and consist of the requirements of (i) stabilizability of the pair of matrices for the equations of the follower objects, (ii) Hurwitz character and (iii) coincidence of the matrices for the leaders in the case of numerous leaders (the multileader case), (iv) solvability of some linear
equations and constraints (of the type of equalities) upon the vectors defining the desired mutual location between the leader and the follower. Furthermore, the whole class of controls, which provide for satisfaction of the linear internal stability, has been described. While relying upon the conditions obtained, it was demonstrated that practically only one-leader formations can possess internal stability. We have identified a subclass (the formations whose graph is the input tree) in the class with one leader, in which there are no constraints of the type of equalities,  that represent the main problem (obstacle) for the internal stability of multi-leader formations.

Keywords: formation, input-to-state stability, acyclic digraph, Hurwitz matrix, sta\-bi\-li\-za\-bi\-li\-ty



\end{document}